\documentclass[a4paper, 10pt, leqno]{amsbook}

\usepackage[T1]{fontenc}
\usepackage[latin1]{inputenc}
\usepackage{ngerman}
\usepackage{amsmath}
\usepackage{amssymb}
\usepackage{setspace}

\newcounter{nr}
\newenvironment{nummer}[2]{\newcounter{#2}\setcounter{#2}{\value{nr}}\textbf{\textsc{#1 \arabic{chapter}.\arabic{nr}.}}\addtocounter{nr}{1}\;\;}{}
\newenvironment{beweis}{\textsc{Proof:}\;\;}{\\ \rightline{$\Box$}\\}

\newenvironment{litfwd}[1]{\newcounter{#1}\setcounter{#1}{\value{nr}}\addtocounter{nr}{1}\;\;}{}

\begin{document}

\setcounter{nr}{1}

\begin{litfwd}{cab}\end{litfwd}
\begin{litfwd}{chsa}\end{litfwd}

\begin{litfwd}{fer}\end{litfwd}

\begin{litfwd}{ham}\end{litfwd}
\begin{litfwd}{hit1}\end{litfwd}
\begin{litfwd}{hit2}\end{litfwd}

\centerline{} \vspace{1cm}

\centerline{\textsc{\Large \textbf{LIFTING $SU(3)$-STRUCTURES}}}
\centerline{\textsc{\Large \textbf{TO NEARLY PARALLEL
$G_{2}$-STRUCTURES}}}

\vspace{3cm} \centerline{Sebastian Stock}\centerline{Mathematical
Institute} \centerline{University of
Cologne}\centerline{stock.s@uni-koeln.de}

\vspace{2cm}

\centerline{\textsc{\textbf{Abstract}}}
\begin{center}
\parbox{10cm}{Hitchin shows in [\arabic{hit1},\arabic{hit2}] that
half-flat $SU(3)$-structures on a $6$-dimensional manifold $M$ can
be lifted to parallel $G_{2}$-structure on the product
$M\times\mathbb{R}$. We show that Hitchin's approach can also be
used to construct nearly parallel $G_{2}$-structures by lifting
so-called nearly half-flat structures. These $SU(3)$-structures
are described by pairs $(\omega,\varphi)$ of stable $2$- and
$3$-forms with $d\varphi=\lambda\omega^{2}$ for some
$\lambda\in\mathbb{R}\setminus\{0\}$.}
\end{center}

\centerline{} \vspace{2cm}
{\large \textbf{Introduction}}\\

The group $SU(3)$ can be realized as the $G_{2}$-stabilizer of a
point in the $6$-sphere. Therefore a $SU(3)$-structure on a
$6$-dimensional manifold $M$ can always be lifted to a
$G_{2}$-structure on the product $M\times\mathbb{R}$. Conversely a
$G_{2}$-structure induces a $SU(3)$-structure on oriented
hypersurfaces (see [\arabic{cab}]). This relationship between
$SU(3)$ and $G_{2}$-structures seems to be very pronounced in the
case of (nearly) half-flat $SU(3)$- and (nearly) parallel
$G_{2}$-structures. Cabrera shows in [\arabic{cab}] that (nearly)
parallel $G_{2}$-structures induce (nearly) half-flat
$SU(3)$-structures on oriented hypersurfaces. Conversely Hitchin
shows in [\arabic{hit1},\arabic{hit2}] that half-flat
$SU(3)$-structures can be lifted to parallel $G_{2}$-structures.
It is well known that the associated metric of a parallel
$G_{2}$-structure is Ricci-flat. Nearly parallel $G_{2}$-structure
can be described by special $3$-forms $\psi$ with
$d\psi=\lambda\ast\psi$ for some
$\lambda\in\mathbb{R}\setminus\{0\}$. The associated metric of a
nearly parallel $G_{2}$-structure is Einstein with constant scalar
curvature. Therefore the question arises (see [\arabic{fer}])
whether nearly half-flat $SU(3)$-structures can
be lifted to nearly parallel $G_{2}$-structures.\\\\

\setcounter{chapter}{1}\setcounter{nr}{1}
\newcounter{chap1}\setcounter{chap1}{\value{chapter}}
{\large \textbf{1. Stable forms and $SU(3)$-structures}}\\

\begin{nummer}{Definition}{defstab}\end{nummer}
Let $V$ be a real $n$-dimensional vector space and
$\rho\in\Lambda^{k}V^{\ast}$. We say that $\rho$ is stable if the
orbit of $\rho$ under the natural action of $GL:=GL(V)$ on
$\Lambda^{k}V^{\ast}$ is open. Consequently we call a $k$-Form
$\rho:FM\longrightarrow\Lambda^{k}\mathbb{R}^{n\ast}=:\Lambda^{k}$
on a $n$-dimensional manifold $M$ with frame bundle $FM$ stable if
for every frame $p\in FM$ the $GL(n,\mathbb{R})$-orbit
$GL.\rho(p)$ of $\rho(p)$ is open in $\Lambda^{k}$. Stable forms
have the property that all forms in a sufficiently small
neighbourhood are stable, too.\\


\begin{nummer}{Example}{exastab}\end{nummer}
Consider the following forms on $\mathbb{R}^{6}$
\begin{enumerate}
\item[] $\omega_{0}:=e^{14}+e^{25}+e^{36}$ \item[]
$\varphi_{0}:=e^{123}-e^{156}+e^{246}-e^{345}$ \item[]
$\sigma_{0}:=\frac{1}{2}\omega_{0}^{2}=e^{1425}+e^{1436}+e^{2536}$
\end{enumerate}
where $(e_{1},..,e_{6})$ denotes the standard basis of
$\mathbb{R}^{6}$. The $2$-form $\omega_{0}$ is the canonical
symplectic form on $\mathbb{R}^{6}$ and its stabilizer under the
$GL(6,\mathbb{R})$ action is $Sp(6,\mathbb{R})$. It follows
$$\text{dim}(GL(6,\mathbb{R}))-\text{dim}(Sp(6,\mathbb{R}))=36-21=\text{dim}(\Lambda^{2}\mathbb{R}^{6\ast})$$
and therefore $\omega_{0}$ is stable. The identity component of
the stabilizer of $\varphi_{0}$, $\sigma_{0}$ is
$SL(3,\mathbb{C})$, $Sp(6,\mathbb{R})$ respectively. Therefore
$\varphi_{0}$ and $\sigma_{0}$ are stable forms too.\\


\begin{nummer}{Remark}{remstab}\end{nummer} The $3$-forms lying in
the open orbit of $\varphi_{0}$ are of special interest. Given
such a form $\varphi$ it is possible (see [\arabic{hit1}]) to
construct a complex structure $I(\varphi)\in
\text{End}(\mathbb{R}^{6})$ on $\mathbb{R}^{6}$. To construct an
almost complex-structure on a $6$-dimensional manifold $M$
consider the subset
$$\Omega_{\varphi_{0}}(M):=\{\varphi:FM\longrightarrow GL.\varphi_{0}\subset\Lambda^{3}\}\subset\Omega^{3}(M)$$
of $3$-forms on $M$. Given $\varphi\in\Omega_{\varphi_{0}}(M)$ we
obtain an almost complex structure
$I(\varphi)\in\Gamma(\text{End}(TM))$. To be more precise, we
require that $M$ is oriented to construct the almost complex
structure. This is due to the fact that only the identity
component of the stabilizer of $\varphi$ is equal to
$SL(3,\mathbb{C})\subset GL(3,\mathbb{C})$. Apart from forms of
type $\varphi_{0}$ there is only one more type of stable $3$-forms
on a $6$-dimensional vector space [\arabic{hit1}]. Henceforth a
stable $3$-form shall always be a form of type $\varphi_{0}$.\\


Let $\rho_{0}\in\{\omega_{0},\varphi_{0},\sigma_{0}\}$ be a stable
$k$-form on $\mathbb{R}^{6}$. In the appendix we give a definition
of an associated volume
$$\epsilon:GL.\rho_{0}\longrightarrow\Lambda^{6}$$
which is compatible with the action of $GL(6,\mathbb{R})$ in the
following way:\\


\begin{nummer}{Lemma}{defvol}\end{nummer}
For $\rho\in GL.\rho_{0}\subset\Lambda^{k}$ and $A\in
GL(6,\mathbb{R})$ we have
\begin{displaymath}
\epsilon(A.\rho)=
\begin{cases}
A.\epsilon(\rho) &\text{ if $\rho\in GL.\omega_{0}$} \\
\text{sgn}(A)A.\epsilon(\rho) &\text{ if $\rho\in GL.\varphi_{0}$ or $\rho\in GL.\sigma_{0}$}\\
\end{cases}
\end{displaymath}
where $\text{sgn}(A)$ denotes the sign of $\text{det}(A)$.\\


Applying this invariance to the action of the scalar matrices, we
see that $\epsilon(\lambda^{k}\rho)=\lambda^{6}\epsilon(\rho)$,
i.e. $\epsilon$ is homogeneous of degree $6/k$ and by Eulers
formula we get $D_{\rho}\epsilon(\rho)=\frac{6}{k}\epsilon(\rho)$.
Finally the wedge product gives an isomorphism $\Lambda^{6-k}\cong
\text{Hom}(\Lambda^{k},\Lambda^{6})$ and for every $\rho\in
GL.\rho_{0}$ we find a unique $\widehat{\rho}\in\Lambda^{6-k}$ for
which
$$D_{\rho}\epsilon(.)=\frac{1}{2}\widehat{\rho}\wedge.$$
holds. By definition of $\widehat{\rho}$ we also have
$$\epsilon(\rho)=\frac{k}{12}\widehat{\rho}\wedge\rho$$
and the map
$\;\widehat{.}:GL.\rho_{0}\longrightarrow\Lambda^{n-k}$
satisfies:\\


\begin{nummer}{Lemma}{lemdach}\end{nummer} For $\rho\in
GL.\rho_{0}$ and $A\in GL(6,\mathbb{R})$ we have
\begin{displaymath}
\widehat{A.\rho}=
\begin{cases}
A.\widehat{\rho} &\text{ if $\rho\in GL.\omega_{0}$. }\\
\text{sgn}(A)A.\widehat{\rho} &\text{ if $\rho\in GL.\varphi_{0}$ or $\rho\in GL.\sigma_{0}$.}\\
\end{cases}
\end{displaymath}\\
\begin{beweis}
For $\dot{\rho}\in T_{A.\rho}(GL.\rho_{0})\cong\Lambda^{k}$ the
invariance of $\epsilon$ yields
$$D_{A.\rho}\epsilon (\dot{\rho})=\text{sgn}(A)\text{det}(A^{-1})D_{\rho}\epsilon(A^{-1}.\dot{\rho})$$
where $\text{sgn}(A)$ occurs only in the case
$\rho_{0}\in\{\varphi_{0},\sigma_{0}\}$. Therefore we get
\begin{displaymath}
\begin{split}
\widehat{A.\rho}\wedge\dot{\rho}&=D_{A.\rho}\epsilon(\dot{\rho})
=\text{sgn}(A)\text{det}(A^{-1})D_{\rho}\epsilon(A^{-1}.\dot{\rho})\\
&=\text{sgn}(A)A.(\widehat{\rho}\wedge A^{-1}.\dot{\rho})
=\text{sgn}(A)A.\widehat{\rho}\wedge \dot{\rho}
\end{split}
\end{displaymath}
\end{beweis}


In the following proposition we describe how the associated forms
$\widehat{\rho}$ of a stable form $\rho\in GL.\rho_{0}$ look like.
In fact the forms $\widehat{\rho}$ are stable too. According to
Lemma \arabic{chap1}.\arabic{lemdach} we only have to compute
$D_{\rho_{0}}\epsilon$ for
$\rho_{0}\in\{\omega_{0},\varphi_{0},\sigma_{0}\}$ and what we get is:\\


\begin{nummer}{Proposition}{satzdach}\end{nummer} The associated
forms $\widehat{\rho}$ are given by:
\begin{enumerate}
\item[(1)] $\widehat{\omega}=-\omega^{2}$. \item[(2)]
$\widehat{\varphi}=-I(\varphi).\varphi$. \item[(3)]
$\widehat{\sigma}=-\omega\;\;$ if $\sigma=\frac{1}{2}\omega^{2}\in
GL.\sigma_{0}$.\\
\end{enumerate}


Now let $M$ be a $6$-dimensional oriented manifold. Using the
invariance of $\epsilon$ and $\;\widehat{.}\;$ we obtain maps
$\epsilon:\Omega_{\rho_{0}}\longrightarrow\Omega^{n}(M)$ and
$\;\widehat{.}\;:\Omega_{\rho_{0}}\longrightarrow\Omega^{n-k}(M)$.
Stable forms can be used to describe $SU(3)$-structures on $M$. In
particular we have (see [\arabic{hit2}]):\\


\begin{nummer}{Proposition}{thmsu}\end{nummer} There is a one-to-one
correspondence between $SU(3)$-structures on $M$ and pairs
$(\omega,\varphi)$ of stable $2$- and $3$-forms on $M$, which
satisfy the following compatibility conditions:
\begin{enumerate}
\item[(I)] $(\omega,\varphi)$ is positive, i.e.
$\omega(X,I(\varphi)X)>0$ for all $X\neq 0$, \item[(II)]
$\omega\wedge\varphi=0$, \item[(III)]
$\epsilon(\varphi)=\epsilon(\omega)$.\\
\end{enumerate}


\begin{nummer}{Remark}{bemsu}\end{nummer} Note that the equation
$\omega\wedge\varphi=0$ holds if and only if
$\omega\wedge\widehat{\varphi}=0$ holds. In particular this is
exactly the condition for $\omega(.,I(\varphi).)$ to be symmetric.
Therefore the conditions (I) and (II) allow us to reconstruct the
metric of the corresponding $SU(3)$-structure. Equation (III) is
the condition for the $(3,0)$-form
$\alpha:=\varphi+i\widehat{\varphi}$ to be of constant length.\\


\begin{nummer}{Definition}{deftypensu}\end{nummer} We say that a
$SU(3)$-structure $(\omega,\varphi)$ is\\
\begin{displaymath}
\begin{split}
\text{ nearly-Kähler }
&:\Leftrightarrow\;\;\;\;\;\;d\omega=3\varphi\;\;\;\text{ and }\;\;\;d\widehat{\varphi}+2\omega^{2}=0,\\
\text{ half-flat }
&:\Leftrightarrow\;\;\;\;\;\;d\varphi=0\;\;\;\text{ and
}\;\;\;d\omega^{2}=0,\\
\text{ nearly half-flat }
&:\Leftrightarrow\;\;\;\;\;\;\exists\;\lambda\in\mathbb{R}\setminus\{0\}:\;\;\; d\varphi+\lambda\omega^{2}=0.\\
\end{split}
\end{displaymath}\\
Every nearly-Kähler structure $(\omega,\varphi)$ is half-flat and
the corresponding $SU(3)$-structure $(\omega,\widehat{\varphi})$
is nearly half-flat. Note that half-flat structures are not
necessarily nearly half-flat. The reason (see [\arabic{chsa}]) to
call structures with $d\varphi=\lambda\omega^{2}$ nearly half-flat
is the type of their intrinsic torsion (in the
Gray-Hervella-decomposition). Apart from the half-flat components
$\mathcal{W}_{1}^{-}\oplus\mathcal{W}_{2}^{-}\oplus\mathcal{W}_{3}$
the torsion of a nearly half-flat structure has additional values
only in $\mathcal{W}_{1}^{+}$, given by
$0\neq\lambda\in\mathbb{R}\cong\mathcal{W}_{1}^{+}$.\\


\newpage
\setcounter{chapter}{2}\setcounter{nr}{1}
\newcounter{chap2}\setcounter{chap2}{\value{chapter}}
{\large \textbf{2. Lifting nearly half-flat $SU(3)$-structures to nearly parallel $G_{2}$-structures}}\\

A $G_{2}$-structure on a $7$-dimensional manifold $M$ can be
described by a stable $3$-form $\psi$ which is locally of the form
$$\psi_{0}:=\omega_{0}\wedge
e^{7}+\varphi_{0}=e^{147}+e^{257}+e^{367}+e^{123}-e^{156}+e^{246}-e^{345}$$
where $(e_{1},..,e_{7})$ denotes the standard basis of
$\mathbb{R}^{7}$. The stabilizer of $\psi_{0}$ under the natural
action of $GL(7,\mathbb{R})$ is the $14$-dimensional Lie-group
$G_{2}$. Choosing a point $x\in S^{6}$ the group $SU(3)$ can be
realized as the subgroup of $G_{2}$ which leaves $x$ invariant.
Due to this fact we have:\\


\begin{nummer}{Lemma}{lemg2}\end{nummer} Suppose $(\omega_{t},\varphi_{t})_{t\in
I}$ is a family of $SU(3)$-structures on a $6$-dimensional
manifold $M^{6}$. Define a $3$-form $\psi$ on $M^{7}:=M^{6}\times
I$ by
$$\psi:=\omega\wedge dt + \varphi,$$
i.e. for $(m,t)\in M^{7}$ we have $\psi_{(m,t)}=\omega_{t,m}\wedge
dt + \varphi_{t,m}$. Then $\psi$ defines a $G_{2}$-structure on
$M^{7}$ and for the induced metric and orientation on $M^{7}$ the
equations
\begin{enumerate}
\item[(1)] $\ast\psi=-\widehat{\varphi}\wedge dt-
\frac{1}{2}\omega^{2}$, \item[(2)] $d^{7}\psi=(d^{6}\omega-
\dot{\varphi})\wedge dt +d^{6}\varphi$, \item[(3)]
$d^{7}\ast\psi=(-d^{6}\widehat{\varphi}-\omega\wedge\dot{\omega})\wedge
dt - \omega\wedge d^{6}\omega$
\end{enumerate}
hold.\\


As a direct consequence of Lemma \arabic{chap2}.\arabic{lemg2} we get:\\


\begin{nummer}{Lemma}{lemnhfnp}\end{nummer} Suppose $(\omega_{t},\varphi_{t})_{t\in
I}$ is a family of $SU(3)$-structures on a 6 dimensional manifold
$M$ and let $\lambda\in\mathbb{R}\setminus\{0\}$. The induced
$G_2$-structure $\psi=\omega\wedge dt+\varphi$ on $M\times I$ is
nearly parallel with $d\psi=\lambda\ast\psi$ if the following
Evolution Equations hold:
\begin{align}
d\varphi_{t}+\frac{\lambda}{2}\omega_{t}^{2}&=0.\\
d\omega_{t}-\dot{\varphi}_{t}&=-\lambda\widehat{\varphi}_{t}.
\end{align}\\


Let $M$ be a closed oriented manifold with a nearly half-flat
$SU(3)$-structure $(\omega,\varphi)$ which satisfies
$$d\varphi=-\frac{\lambda}{2}\omega^{2}$$
for some $\lambda\in\mathbb{R}\setminus\{0\}$. Consider the
non-empty set
$$\mathcal{A}:=\{\varphi\in\Omega^{3}(M)\mid
\exists\;\omega\in\Omega_{\omega_{0}}(M):\;
 d\varphi=-\frac{\lambda}{2}\omega^{2}\}.$$
As a first order differential operator the exterior derivative
$d:\Omega^{4}(M)\longrightarrow\Omega^{3}(M)$ is a continuous map
and the set
$\mathcal{A}=(-\frac{1}{\lambda}d)^{-1}(\Omega_{\sigma_{0}}(M))\subset
\Omega^{3}(M)$ is open in the Fr\'{e}chet topology (see
[\arabic{ham}]). From $\mathcal{A}\neq\emptyset$ we get
$$T_{\varphi}\mathcal{A}=\Omega^{3}(M).$$\\


\begin{nummer}{Lemma}{defpi}\end{nummer}
According to Proposition \arabic{chap1}.\arabic{satzdach} the map
$$\pi:\mathcal{A}\longrightarrow\Omega_{\omega_{0}}(M)\;\;\;\text{ with }\;\;\;\varphi\longmapsto -\;\;\widehat{(-\frac{1}{\lambda}d\varphi)}$$
has the following properties:
\begin{enumerate}
\item[(i)] $-\frac{\lambda}{2}\pi(\varphi)^{2}=d\varphi$.
\item[(ii)] For $\dot{\varphi}\in
T_{\varphi}\mathcal{A}=\Omega^{3}(M)$ we have
$\pi(\varphi)\wedge\pi_{\ast\varphi}(\dot{\varphi})=-\frac{1}{\lambda}d\dot{\varphi}$.
\item[(iii)] For $F\in \text{Diff}(M)$ we have
$\pi(F^{\ast}\varphi)=F^{\ast}\pi(\varphi)$.\\
\end{enumerate}


\begin{nummer}{Definition}{defhom}\end{nummer} Suppose
$(\omega,\varphi)$ is a nearly half-flat $SU(3)$-structure on a
closed manifold $M$ with $d\varphi=-\frac{\lambda}{2}\omega^{2}$
for some $\lambda\in\mathbb{R}\setminus\{0\}$. On the open set
$$A:=\mathcal{A}\cap\Omega_{\varphi_{0}}(M)\subset\Omega^{3}(M)$$
we define a non-degenerated skew-symmetric bilinear form $\Omega$
by
$$\Omega(\dot{\varphi}_{1},\dot{\varphi}_{2}):=\int_{M}\dot{\varphi}_{1}\wedge\dot{\varphi}_{2}.$$
Consider the real-valued function $H:A\longrightarrow \mathbb{R}$
defined by
$$\varphi\longmapsto 2\lambda\left(\int_{M}\epsilon(\varphi)-\int_{M}\epsilon(\pi(\varphi))\right).$$
From Lemma 2.\arabic{defpi} (iii) and the invariance of $\epsilon$
we get
$$H(F^{\ast}\varphi)=\pm H(\varphi)$$
for each $F\in \text{Diff}(M)$. The sign is positive if $F$ is
orientation-preserving, and negative if $F$ changes the
orientation. Now let $X$ be the associated vector field on $A$,
i.e.
$$\Omega(X,.)=DH$$
and let $I\subset\mathbb{R}$ be an open interval for which the
flow $\Phi$ of the vector field $X$ through $\varphi\in A$ is
defined for all times $t\in I$. Then for all $t\in I$ the forms
$$\varphi_{t}:=\Phi_{t}(\varphi)\;\;\;\text{ and }\;\;\;\omega_{t}:=\pi(\varphi_{t})$$
are stable and our main result is:


\begin{nummer}{Theorem}{thmhitnhf}\end{nummer} For a sufficiently
small interval $I$ the forms $(\omega_{t},\varphi_{t})_{t\in I}$
define a family of $SU(3)$-structures on $M$ and the induced
$G_{2}$-structure $\psi=\omega\wedge dt+ \varphi$ on $M\times I$
is nearly parallel with $d\psi=\lambda\ast\psi$.\\\\
\begin{beweis}
The forms $(\omega_{t},\varphi_{t})$ are stable and according to
Proposition 1.\arabic{thmsu} the forms $(\omega_{t},\varphi_{t})$
define a $SU(3)$-structure if the conditions
$$\omega_{t}\wedge\varphi_{t}=0\;\;\;\text{ and }\;\;\;\epsilon(\omega_{t})=\epsilon(\varphi_{t})$$
are satisfied. Note that the open condition (I) of Proposition
1.\arabic{thmsu} holds on the compact manifold $M$ for sufficient
small $I$ and Lemma 2.\arabic{defpi} (i) yields
\setcounter{equation}{0}
\begin{equation}
d\varphi_{t}+\frac{\lambda}{2}\omega_{t}^{2}=0.
\end{equation}
Therefore the induced $G_{2}$-structure is nearly parallel if the
equation
\begin{equation}
d\omega_{t}-\dot{\varphi}_{t}=-\lambda\widehat{\varphi}_{t}
\end{equation}
holds (see Lemma 2.\arabic{lemnhfnp}). For $Y\in\Gamma(TM)$
consider the map
$$\mu_{Y}:A\longrightarrow\mathbb{R}\;\;\;\text{ with }\;\;\;
\varphi\longmapsto\int_{M}(Y\lrcorner\;(-\frac{1}{\lambda}d\varphi))\wedge\varphi.$$
For $\varphi\in A$ we have $d\varphi=-\frac{\lambda}{2}\omega^{2}$
and therefore
$$\mu_{Y}(\varphi)=\int_{M}(Y\lrcorner\;\omega)\wedge\omega\wedge\varphi.$$
Since $\omega$ is non-degenerated, the equation
$\omega_{t}\wedge\varphi_{t}=0$ is satisfied if and only if
$\mu_{Y}(\varphi_{t})=0$ holds for all $Y\in\Gamma(TM)$, i.e.
\begin{equation}
d\mu_{Y}(\dot{\varphi}_{t})=0
\end{equation}
for all $Y\in\Gamma(TM)$. Stokes' theorem yields
\begin{displaymath}
\begin{split}
\Omega(X\circ\varphi_{t},L_{Y}\varphi_{t})&=\int_{M}\dot{\varphi}_{t}\wedge L_{Y}\varphi_{t}\\
&=\int_{M}\dot{\varphi}_{t}\wedge
d(Y\lrcorner\varphi_{t})+\dot{\varphi}_{t}\wedge(Y\lrcorner
d\varphi_{t})\\
&=\int_{M}d\dot{\varphi}_{t}\wedge (Y\lrcorner\varphi_{t})
+\dot{\varphi}_{t}\wedge(Y\lrcorner d\varphi_{t}).
\end{split}
\end{displaymath}
and from $0=d\dot{\varphi}_{t}\wedge\varphi_{t}$ we get
\begin{displaymath}
\begin{split}
\Omega(X\circ\varphi_{t},L_{Y}\varphi_{t})&=\int_{M}-(Y\lrcorner
d\dot{\varphi}_{t})\wedge\varphi_{t}+\dot{\varphi}_{t}\wedge(Y\lrcorner d\varphi_{t})\\
&=\int_{M}-(Y\lrcorner d\dot{\varphi}_{t})\wedge\varphi_{t}
-(Y\lrcorner d\varphi_{t})\wedge\dot{\varphi}_{t}\\
&=\lambda d\mu_{Y}(\dot{\varphi}_{t}).
\end{split}
\end{displaymath}
Hence
$$\lambda
d\mu_{Y}(\dot{\varphi}_{t})=\Omega(X\circ\varphi_{t},L_{Y}\varphi_{t})=D_{\varphi_{t}}H(L_{Y}\varphi_{t}).$$\\
On the compact manifold $M$ let $\Psi_{s}$ be the flow of $Y$,
defined for $|s|<\varepsilon$. For $\varphi_{t}\in A$ define a
path
$$c:(-\varepsilon,\varepsilon)\longrightarrow A\;\;\;\text{ by
}\;\;\;s\longmapsto\Psi_{s}^{\ast}\varphi_{t}.$$\\
Note that
$dc(s)=\Psi_{s}^{\ast}d\varphi_{t}=-\frac{\lambda}{2}(\Psi_{s}^{\ast}\omega_{t})^{2}$
holds. Therefore we have $c(s)\in A$ with
$\dot{c}(0)=L_{Y}\varphi_{t}$ for sufficiently small
$\varepsilon$.
Finally, Definition 2.\arabic{defhom} yields (3):\\
\begin{displaymath}
\begin{split}
\lambda d\mu_{Y}(\dot{\varphi}_{t})&=D_{\varphi_{t}}H(L_{Y}\varphi_{t})\\
&=(s\longmapsto H\circ c(s))'(0)\\
&=(s\longmapsto \underbrace{H(\Psi_{s}^{\ast}\varphi_{t})}_{\equiv
H(\varphi_{t})})'(0)=0.
\end{split}
\end{displaymath}
So we have $\omega_{t}\wedge\varphi_{t}=0$ for all $t\in I$ or
equivalent
\begin{equation}
\omega_{t}\wedge\widehat{\varphi}_{t}=0.
\end{equation}
By definition we have
$\epsilon(\pi(\varphi))=-\frac{1}{6}\pi(\varphi)^{3}$ and the flow
equation $\Omega(X,.)=DH$ yields
\begin{displaymath}
\begin{split}
\int_{M}\dot{\varphi}_{t}\wedge\dot{\varphi}&=\Omega(\dot{\varphi}_{t},\dot{\varphi})
=\Omega(X\circ\varphi_{t},\dot{\varphi})=D_{\varphi_{t}}H(\dot{\varphi})\\
&=2\lambda\int_{M}D_{\varphi_{t}}\epsilon(\dot{\varphi})+\frac{1}{2}\pi(\varphi_{t})^{2}\wedge\pi_{\ast}(\dot{\varphi})
\end{split}
\end{displaymath}
for all $\dot{\varphi}\in T_{\varphi_{t}}A=\Omega^{3}(M)$. Using
Lemma 2.\arabic{defpi} (ii) we obtain Equation (2):
\begin{displaymath}
\begin{split}
\int_{M}\dot{\varphi}_{t}\wedge\dot{\varphi}&=2\lambda\int_{M}\frac{1}{2}\widehat{\varphi}_{t}\wedge\dot{\varphi}
+\frac{1}{2}\omega_{t}\wedge(-\frac{1}{\lambda}d\dot{\varphi})\\
&=\int_{M}\lambda\widehat{\varphi}_{t}\wedge\dot{\varphi}-\omega_{t}\wedge
d\dot{\varphi}\\
&=\int_{M}(\lambda\widehat{\varphi}_{t}+d\omega_{t})\wedge
\dot{\varphi}.\\
\end{split}
\end{displaymath}
Further we have
\begin{displaymath}
\begin{split}
2(D\epsilon(\dot{\omega}_{t})-D\epsilon(\dot{\varphi}_{t}))&=
\widehat{\omega}_{t}\wedge\dot{\omega}_{t}-\widehat{\varphi}_{t}\wedge\dot{\varphi}_{t}\\
&\overset{(2)}{=}-\omega_{t}^{2}\wedge\dot{\omega}_{t}-\widehat{\varphi}_{t}\wedge(\lambda\widehat{\varphi}_{t}+d\omega_{t})\\
&=-\omega_{t}^{2}\wedge\dot{\omega}_{t}-\widehat{\varphi}_{t}\wedge
d\omega_{t}\\
&\overset{(4)}{=}-\omega_{t}^{2}\wedge\dot{\omega}_{t}-d\widehat{\varphi}_{t}\wedge\omega_{t}\\
&\overset{(2)}{=}-\omega_{t}^{2}\wedge\dot{\omega}_{t}-\frac{1}{\lambda}d\dot{\varphi}_{t}\wedge\omega_{t}\\
&\overset{(1)}{=}0
\end{split}
\end{displaymath}
which proves the theorem.\\
\end{beweis}


According to the results in [\arabic{cab}] and
[\arabic{hit1},\arabic{hit2}] Theorem
\arabic{chap2}.\arabic{thmhitnhf} is the last part to completely
understand the relationship between (nearly) half-flat
$SU(3)$-structures on $M$ and (nearly) parallel $G_{2}$-structures
on the product $M\times\mathbb{R}$.\\
Apart from nearly half-flat structures on hypersurfaces of nearly
parallel $G_{2}$ manifolds, nearly Kähler structures are also
examples for nearly half-flat structures (see [\arabic{fer}]). A
nearly Kähler structure is a $SU(3)$-structure $(\omega,\varphi)$
which satisfies
$$d\omega=3\varphi\;\;\;\text{ and }\;\;\;d\widehat{\varphi}+2\omega^{2}=0.$$
These equations imply that $(\omega,\varphi)$ is a half-flat
structure and the corresponding $SU(3)$-structure
$(\omega,\widehat{\varphi})$ is nearly half-flat. In the
nearly-Kähler case Hitchin's evolution Equations and the evolution
Equations (1) and (2) of Lemma \arabic{chap2}.\arabic{lemnhfnp}
can be solved explicitly in terms of $(\omega,\varphi)$ and
$(\omega,\widehat{\varphi})$: Suppose $(\omega,\varphi)$ is a
nearly-Kähler structure with metric $g$. For $t\in\mathbb{R}_{>0}$
the $SU(3)$-family
$$\omega_{t}:=t^{2}\omega\;\;\;\text{ and }\;\;\;\varphi_{t}:=t^{3}\varphi$$
induces a parallel $G_{2}$-structure on the product
$M\times\mathbb{R}_{>0}$ and the $G_{2}$-metric $\widehat{g}$ is
given by the conical metric $\widehat{g}=t^{2}g+dt^{2}$ on
$M\times\mathbb{R}_{>0}$. We can do something similar to construct
nearly parallel $G_{2}$-structures: For $t\in(0,\pi)$ the
$SU(3)$-family
$$\omega_{t}:=\sin^{2}(t)\omega\;\;\;\text{ and
}\;\;\;\varphi_{t}:=\sin^{3}(t)(\sin(t)\widehat{\varphi}+\cos(t)\varphi)$$
induces a nearly parallel $G_{2}$-structure on the product
$M\times(0,\pi)$ and the $G_{2}$-metric $\widehat{g}$ is given by
$\widehat{g}=\sin^{2}(t)g+dt^{2}$.\\

Hitchin's construction and the construction presented in Chapter
$2$ have in common that the torsion type of the $SU(3)$-family
remains constant. This indicates that a unifying approach to both
constructions might be possible.\\\\\\

{\large \textbf{Appendix: Definition of the volumes} (see [\arabic{hit2}])}\\

For the case n=6 and k=2 the volume associated to $\omega\in
GL.\omega_{0}$ is simply the Liouville volume
$$\epsilon(\omega):=-\frac{1}{3!}\omega^{3}.$$\\
In the other cases, we construct $GL(6,\mathbb{R})$-invariant maps
in the following way:\\

\textbf{(i) n=6, k=3:}
$$K:\Lambda^{3}\longrightarrow
\text{Hom}(\mathbb{R}^{6},\mathbb{R}^{6}\otimes\Lambda^{6})\;\;\;\text{
by }\;\;\; K_{\varphi}(x)(\alpha):=\frac{1}{2}\alpha\wedge
(x\lrcorner\varphi)\wedge\varphi\in \Lambda^{6}$$
($x\in\mathbb{R}^{6}$, $\alpha\in\mathbb{R}^{6\ast}$) and
$$\lambda:\Lambda^{3}\longrightarrow (\Lambda^{6})^{2}\;\;\; \text{
by } \;\;\;\varphi\longmapsto
\frac{1}{6}\text{tr}(K^{2}_{\varphi}).$$\\
\textbf{(ii) n=6, k=4:}
$$K:\Lambda^{4}\longrightarrow
\text{Hom}(\mathbb{R}^{6\ast},\mathbb{R}^{6}\otimes\Lambda^{6})\;\;\;
\text{ by }
\;\;\;K_{\sigma}(\alpha)(\beta):=\beta\wedge\alpha\wedge\sigma\in
\Lambda^{6}$$ ($\alpha,\beta\in\mathbb{R}^{6\ast}$) and
$$\lambda:\Lambda^{4}\longrightarrow (\Lambda^{6})^{4}\;\;\; \text{ by }
\;\;\;\sigma\longmapsto \text{det}(K_{\sigma}).$$\\

Note that $\lambda(\varphi_{0})<0$ and $\lambda(\sigma_{0})>0$
holds. Therefore we obtain the associated volumes by:

\begin{enumerate}
\item[(i)]$\epsilon:GL.\varphi_{0}\longrightarrow\Lambda^{6}\;\;\;\text{
with }\;\;\;\varphi\longmapsto (-\lambda(\varphi))^{\frac{1}{2}}.$
\item[(ii)]$\epsilon:GL.\sigma_{0}\longrightarrow\Lambda^{6}\;\;\;\text{
with }\;\;\;\sigma\longmapsto (\lambda(\sigma))^{\frac{1}{4}}.$
\end{enumerate}

The complex structure associated to $\varphi\in GL.\varphi_{0}$ is
given by
$$I(\varphi):=\frac{1}{\epsilon(\varphi)}K(\varphi)\in \text{End}(\mathbb{R}^{6}).$$

\newpage\setcounter{nr}{1}\pagestyle{empty}
\textsc{\large \textbf{References}}\\


\textsc{\textbf{[\arabic{cab}]} M. Cabrera:
}\textit{$SU(3)$-structures on hypersurfaces of manifolds with
$G_{2}$-structure}, 2005, to appear in Monatsh. Math. Preprint
math.DG/0410610 v6.\\

\textsc{\textbf{[\arabic{chsa}]} S.Chiossi, S.Salamon:
}\textit{The intrinsic torsion of $SU(3)$ and $G_{2}$ structures},
Differential Geometry, Valencia 2001, World Sci. Publishing, 2002,
pp 115-133, math.DG/0202282 v1.\\


\textsc{\textbf{[\arabic{fer}]} M. Fernandez, S.Ivanov, V.Munoz,
L.Ugarte: }\textit{Nearly hypo structures and compact nearly
Kähler 6-manifolds with conical singularities}, 2006,
math.DG/ 0602160 v2.\\


\textsc{\textbf{[\arabic{ham}]} R.S. Hamilton: }\textit{The
Inverse Function Theorem of Nash and Moser}, American Math. Soc.,
1982, pages 65-222.\\

\textsc{\textbf{[\arabic{hit1}]} N. Hitchin: }\textit{The geometry
of three-forms in six and seven dimensions}, J. Differ. Geom.,
55:547-576, 2000, math.DG/0010054 v1.\\

\textsc{\textbf{[\arabic{hit2}]} N. Hitchin: }\textit{Stable forms
and special metrics}, Global Differential Geometry: \textit{The
Mathematical Legacy of Alfred Gray}, volume 288 of Contemp. Math.,
pages 70-89. American Math. Soc., 2001, math.DG/0107101 v1.\\

\end{document}